\documentclass[11pt]{amsart}

\newcommand{\area}{{\rm Area}}

\newcommand{\D}{{\mathbb D}}
\newcommand{\id}{{\rm id}}

\usepackage{amsthm,amsmath,amsfonts, amssymb}
\usepackage[all]{xy}
\usepackage{psfrag}
\usepackage[dvips]{graphicx}

\usepackage[latin1]{inputenc}

\makeatletter
\def\th@alexnormal{%
\let\thm@indent\noindent 
\thm@headfont{\bfseries}
\normalfont
}
\makeatother

\makeatletter
\def\th@alexit{%
\let\thm@indent\noindent 
\thm@headfont{\bfseries}
\normalfont
\fontshape{it}
\selectfont
}
\makeatother

\theoremstyle{alexit}

\newtheorem{theorem}[equation]{Theorem}
\newtheorem{proposition}[equation]{Proposition}
\newtheorem{lemma}[equation]{Lemma}

\newtheorem{question}[equation]{Question}

\theoremstyle{remark}
\newtheorem{remark}[equation]{Remark}
\theoremstyle{definition}

\numberwithin{equation}{subsection}

\begin{document}
\author{Alexandre Girouard}
\address{School of Mathematics, Cardiff University,
Senghennydd Road, Cardiff, Wales,  CF24 4AG,
UK}
\email{GirouardA@cardiff.ac.uk}

\author{Iosif Polterovich}\thanks{The second author  is  supported by
NSERC and FQRNT}
\address{D\'epartement de math\'ematiques et de
statistique, Universit\'e de Montr\'eal, C. P. 6128,
Succ. Centre-ville, Montr\'eal, Qu\'ebec,  H3C 3J7, Canada}
\email{iossif@dms.umontreal.ca}

\title[On the Hersch--Payne--Schiffer inequalities]{On the Hersch--Payne--Schiffer inequalities  for
Steklov eigenvalues 
}
\begin{abstract}
We prove that the  isoperimetric inequality due to
Hersch--Payne--Schiffer for the $n$-th nonzero Steklov eigenvalue of a
bounded simply--connected planar domain is sharp for all $n \ge
1$. The equality is attained in the limit by a sequence of
simply--connected domains degenerating to the disjoint union of $n$
identical disks.  We give a new proof of  this  inequality for $n=2$
and show that it is strict in this case.  Related results are also
obtained for the product  of two consecutive Steklov eigenvalues.

\end{abstract}
\maketitle

\section{Introduction and main results}

\subsection{Steklov eigenvalue problem}
Let $\Omega$ be a simply-connected bounded planar domain with
Lipschitz boundary and $\rho\in L^{\infty}(\partial \Omega)$ be a
non-negative nonzero function. The {\it Steklov eigenvalue
problem} \cite{St}  is given by
\begin{equation}
\label{steklov}
\begin{cases}
  \Delta u=0& \mbox{ in } \Omega,\\
  \frac{\partial}{\partial\nu}u=\sigma \,\rho\,u& \mbox{ on
  }\partial\Omega,
\end{cases}
\end{equation}
where $\frac{\partial}{\partial\nu}$ is the outward normal
derivative. There are several physical interpretations of
the Steklov problem \cite{Bandle, Payne}. In particular, it describes the vibration
of a free membrane with its whole mass  $M(\Omega)$ distributed on the boundary with
density $\rho$:
\begin{equation}
\label{mass}
M(\Omega)=\int_{\partial \Omega} \rho(s) ds.
\end{equation}
If $\rho \equiv 1$,   the mass of $\Omega$ is equal to the length of $\partial \Omega$.

The spectrum of the Steklov problem is discrete, and the eigenvalues
$$0=\sigma_0<\sigma_1(\Omega)\leq \sigma_2(\Omega)\leq \sigma_3(\Omega)\leq\cdots\nearrow\infty$$
satisfy the following variational characterization  \cite[p.~95
and p.~103]{Bandle}:
\begin{gather}
\label{rayleigh}
 \sigma_n(\Omega)=\inf_{E_n} \sup_{0\neq u\in E_n}
  \frac{\int_{\Omega}|\nabla u|^2\,dz}
  {\int_{\partial\Omega}u^2\rho\, ds}, \,\,\,\,\, n=1,2,\dots
\end{gather}
Here the infimum is taken over all $n$-dimensional subspaces $E_n$
of the Sobolev space $H^1(\Omega)$ that are orthogonal to constants
on $\partial\Omega$, i.e. $\int_{\partial\Omega}u(s)\rho(s)\,ds=0$
for all $u\in E_n$.
 Note that, as in the case of Neumann boundary conditions, the Steklov spectrum always starts with the eigenvalue $\sigma_0=0$,
and the corresponding eigenfunctions are constant.

If  the density $\rho\equiv 1$,  the Steklov eigenvalues and
eigenfunctions  coincide with those of the
\emph{Dirichlet-to-Neumann operator} $$\Gamma:
H^{1/2}(\partial\Omega)\rightarrow H^{-1/2}(\partial\Omega)$$
defined by
$$\Gamma f=\frac{\partial}{\partial\nu}(\mathcal{H}f),$$
where $\mathcal{H}f$ is the unique harmonic extension of the
function $f \in H^{1/2}(\partial\Omega)$ to the interior of
$\Omega$. If the boundary is smooth, the Dirichlet-to-Neumann
operator is a first order elliptic pseudo-differential operator
\cite[pp. 37--38]{Taylor}.  It has  various important applications,
particularly to the study of inverse problems \cite{US}.

\subsection{Upper bounds on Steklov eigenvalues}
The present paper is motivated by the following
\begin{question}
\label{quest}
How large can the $n$-th eigenvalue of the Steklov problem be on a bounded simply-connected planar domain of a given mass?
\end{question}

\smallskip

For $n=1$, the answer to Question \ref{quest}  was given in 1954 by Weinstock  \cite{Weinst}.  He proved that
\begin{equation}
\label{Weinst} \sigma_1(\Omega)\, M(\Omega) \leq  2\pi
\end{equation}
with the equality attained on a disk with $\rho \equiv {\rm const}$.
Note that the first eigenvalue of the unit disk $\D$ with $\rho
\equiv 1$ has multiplicity two and $\sigma_1(\D)=\sigma_2(\D)=1$.
Various extensions of Weinstock's inequality and related results
can be found in  \cite{Bandle2, HPS2, Brock, Dittmar, Edward}; see
also \cite[section~8]{AB} for a recent survey.

\smallskip

In 1974, Hersch--Payne--Schiffer \cite[p. 102]{HPS1} proved the
following estimates:
\begin{gather}\label{hpsGeneral}
  \sigma_m(\Omega)\,\sigma_n(\Omega)\,M(\Omega)^2\leq
  \begin{cases}
    (m+n-1)^2\,\pi^2&\mbox{ if } m+n \mbox{ is odd},\\
    (m+n)^2\,\pi^2&\mbox{ if } m+n \mbox{ is even}.
  \end{cases}
\end{gather}
In particular, for $m=n$ and $m=n+1$ we get 
\begin{equation}
\label{hps} \sigma_n(\Omega)M(\Omega)\,  \le 2 \pi n, \,\,\,
n=1,2,\dots,
\end{equation}
\begin{equation}
\label{hps2}
\sigma_n(\Omega)\,\sigma_{n+1}(\Omega)  M(\Omega)^2 \le 4\pi^2 n^2, \,\,\, n=1,2,\dots.
\end{equation}
 \subsection{Main results}  If  $n=1$,  it is easy to see that~\eqref{hps}  and \eqref{hps2}
become equalities on a disk with constant density $\rho$ on the boundary. It was remarked in \cite{HPS1} that  estimates 
\eqref{hpsGeneral} are not expected to be sharp for all  $m$ and $n$. While this is likely to be true, it turns out that if  $m=n$ or  $m=n+1$,  these inequalities  {\it are} sharp for all $n\ge 1$.
\begin{theorem}
  \label{maintheorem}
  There  exists a family of simply-connected bounded
  Lipschitz domains $\Omega_\varepsilon  \subset \mathbb{R}^2$,
  with $\rho \equiv 1$ on $\partial \Omega_\varepsilon$ for all
  $\varepsilon$, degenerating to the disjoint union of $n$
  identical disks as $\varepsilon\to 0+$, such that
  \begin{equation}
    \label{limitcase}
    \lim_{\varepsilon\rightarrow 0+}\sigma_n(\Omega_\varepsilon) \, M(\Omega_\varepsilon) = 2
    \pi n, \,\,\, n=2,3,\dots
  \end{equation}
  and
  \begin{equation}
    \label{limitcasenew}
    \lim_{\varepsilon\rightarrow 0+}\sigma_n(\Omega_\varepsilon) \, \sigma_{n+1}(\Omega_\varepsilon)\,
    M(\Omega_\varepsilon)^2  = 4
    \pi^2 n^2,  \,\,\, n=2,3,\dots
  \end{equation}
  In particular, the Hersch--Payne--Schiffer inequalities  \eqref{hps}
  and \eqref{hps2} are  sharp for all $n\ge 1$.
\end{theorem}

\begin{remark}
  As we show in subsection~\ref{counterexample}, in order to obtain
  \eqref{limitcase} and \eqref{limitcasenew}, one has to be
  careful in the choice of a family of domains degenerating
  to the disjoint union of $n$ identical disks.
  
  It would be interesting to check whether each of the equalities  \eqref{limitcase} and
  \eqref{limitcasenew}
  {\it implies} that the family $\Omega_\varepsilon$ converges in an
  appropriate sense to the disjoint
  union of $n$ identical disks.
\end{remark}
\begin{remark}
  If  $\rho\equiv 1$,  estimate \eqref{hps}  and  the standard
  isoperimetric inequality in $\mathbb{R}^2$ imply
  $$
  \sigma_n(\Omega)\, \sqrt{\area(\Omega)} < n \sqrt{\pi}, \quad n\ge 2.
  $$
  The sharp ``isoareal'' estimate on $\sigma_n$, $n\ge 2$, is unknown (see \cite[Open Problem~25]{Henrot}).
\end{remark}

Theorem \ref{maintheorem} gives an almost  complete answer  to
Question~\ref{quest}. It remains to establish whether the
inequality~\eqref{hps} is {\it strict} for all $n\ge 2$. We believe it is true. 
A modification of the method introduced in \cite{GNP}  allows to prove this result
for $n=2$.
\begin{theorem}\label{maintheoremII}
  Inequality \eqref{hps} is strict for $n=2$:
  \begin{equation}
    \label{bound:main}
    \sigma_2(\Omega) \, M(\Omega) < 4 \pi.
  \end{equation}
\end{theorem}
The proof of Theorem~\ref{maintheoremII} uses the Riemann mapping
theorem similarly to  \cite{Szego1, Weinst}.  Note that
this  approach is very different from the techniques of
\cite{HPS1}. 

\subsection{Discussion}
To put the inequalities \eqref{Weinst}  and \eqref{bound:main}  in
perspective, let us state similar results for  Dirichlet and
Neumann eigenvalues. Since  Dirichlet and Neumann eigenvalue
problems describe vibrations of a membrane of unit density, the
mass of the membrane is equal to its area.

Let  $\Omega \subset \mathbb{R}^2$ be a bounded domain, and let
$0<\lambda_1(\Omega) \le
\lambda_2(\Omega) \le \dots $  and
$0=\mu_0 <\mu_1(\Omega) \le\mu_2(\Omega)\le \dots$ be the Dirichlet and Neumann
eigenvalues of $\Omega$, respectively. We have:
\begin{itemize}

\item[--]  {\it Faber--Krahn inequality}: $\lambda_1(\Omega) \area(\Omega) \ge \pi\, \lambda_1(\D)$ (conjectured in \cite{Ra},
proved in \cite{Faber} and  \cite{Krahn1}, a weaker version obtained  in  \cite{Co}).

\smallskip

\item[--] {\it Krahn's inequality}:  $\lambda_2(\Omega)
\area(\Omega)> 2\pi \lambda_1(\D)$ (\cite{Krahn2}; see
\cite[p.~110]{AB} for an interesting discussion on the history of
this result). The equality is attained in the limit by a sequence of
domains degenerating to the disjoint union of two identical disks.

\smallskip

\item[--]  {\it Szeg\"o--Weinberger inequality}: $\mu_1(\Omega)
  \area(\Omega)  \le \pi\, \mu_1(\D) $. This estimate  was proved in
  \cite {Szego1} for simply-connected planar domains. In \cite{Wein},
  the result was extended to arbitrary domains in all dimensions.

\smallskip

\item[--] If $\Omega$ is simply-connected, $\mu_2(\Omega)
\area(\Omega)  \le 2 \pi\, \mu_1(\D).$ This inequality  was
recently proved in \cite{GNP}. It is an open question whether it
holds for multiply connected planar domains. The equality is
attained in the limit by a sequence of domains degenerating to the
disjoint union of two identical disks.

\end{itemize}

For higher Dirichlet and Neumann eigenvalues, 
no sharp estimates of this type are known, and the situation is quite different from the Steklov case.
As was mentioned in~\cite[Remark 1.2.8]{GNP},  the disjoint
union of $n$ identical disks can not maximize $\mu_n(\Omega)\,
\area(\Omega)$ for all $n\ge~1$, because this  would contradict Weyl's law. By the same
argument, $\lambda_n(\Omega)\, \area(\Omega)$ is not minimized by
the disjoint union of $n$ disks for $n$ large enough. In fact, it is
conjectured that for $n=3$ the minimizer is a single disk, see
\cite{WK, BH}.
\subsection{Plan of the paper} In Section~\ref{hpssharp} we prove
Theorem~\ref{maintheorem}. We also construct a family of
domains whose Steklov spectrum completely
``collapses'' to zero in the limit as the domains degenerate to the
disjoint union of two unit disks. This phenomenon is quite
surprising and does not occur for either Dirichlet or Neumann
eigenvalues. The rest of the paper is devoted to the proof of Theorem
\ref{maintheoremII}.
In Section \ref{folding:section}, the ``folding and rearrangment''  technique,
introduced in \cite{Nad} and developed in \cite{GNP}, is adapted to
the Steklov problem.  In Section \ref{testfunctions} we combine
analytic and topological arguments to construct a two-dimensional
space of test functions for the variational
characterization~\eqref{rayleigh} of the second Steklov
eigenvalue. This space of test functions is then used to prove inequality
\eqref{bound:main}.


\section{Maximization and collapse of Steklov eigenvalues}
\label{hpssharp}


\subsection{Proof of  Theorem \ref{maintheorem}}
\label{whysharp}
Let us start with the case $n=2$. For each
$\varepsilon \in (0, \frac{1}{10})$,
consider the simply--connected planar domain
\begin{equation}
  \label{family}
  \Omega_\varepsilon=\{|z-1+\varepsilon|<1\}\cup \{|z+1-\varepsilon|<1\} \subset
  {\mathbb C}.
\end{equation}
\begin{figure}[h]
  \centering
  \psfrag{e}[][][1]{$2\varepsilon$}
  \includegraphics[width=11cm]{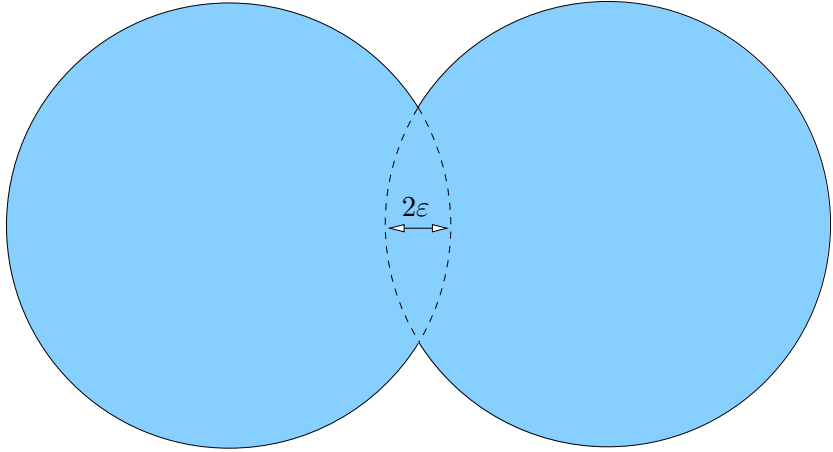}
  \caption{The domain $\Omega_\varepsilon$ for $n=2$}
  \label{figpullappart}
\end{figure}
As $\varepsilon \to 0+$, $\Omega_\varepsilon$ degenerates to the disjoint union of two identical unit
disks.
\begin{lemma}
\label{conv} Let $\rho \equiv 1$ on $\partial \Omega_\varepsilon$ for any $\varepsilon$. Then
\begin{equation*}
\lim_{\varepsilon \to 0+} \sigma_2(\Omega_\varepsilon)=1.
\end{equation*}
\end{lemma}
We recall that if $\rho \equiv 1$ then
$\sigma_1(\D)=\sigma_2(\D)=1$.
\begin{remark}
\label{diff} While this lemma is not surprising, it does not follow
in a straightforward way from general results on convergence of
eigenvalues. The difficulty is that the family $\Omega_\varepsilon$ is
not uniformly Lipschitz. Equivalently, the family $\Omega_\varepsilon$ does not satisfy the
uniform cone condition (see \cite[p.~49]{DZ} or \cite[p.~53]{HP}). This means that one can
not choose the Lipschitz constant uniformly in {\it both} $z\in
\partial \Omega_\varepsilon$ and $\varepsilon$.  Indeed, it is easy
to see that the Lipschitz constant blows up near $z=0$ as
$\varepsilon \to 0$.  In this situation, the Steklov eigenvalues may
apriori have a rather surprising limiting behavior, see subsection
\ref{counterexample}.
\end{remark}

\begin{proof}[Proof of Lemma \ref{conv}] For each $\varepsilon \in (0,\frac{1}{10})$,
$$\sigma_2(\Omega_\varepsilon) M(\Omega_\varepsilon) \le  4\pi $$
by \eqref{hps}. Since  $\lim_{\varepsilon \to 0+}
M(\Omega_\varepsilon)=4\pi$,  we have
$$\limsup_{\varepsilon \to 0+} \sigma_2(\Omega_\varepsilon) \le
1.$$ It remains to show that
\begin{equation}
\label{opposite} \liminf_{\varepsilon \to 0+}
\sigma_2(\Omega_\varepsilon) \ge 1.
\end{equation}

In view of Remark \ref{diff}, in order to apply standard results on
convergence of eigenvalues,  we need to ``desingularize'' the family of
domains $\Omega_\varepsilon$. Let
$\Omega'_\varepsilon=\Omega_\varepsilon \cap\{\Re z < 0 \}$.
Consider the following auxiliary mixed eigenvalue problem on $\Omega'_\varepsilon$: impose
the Neumann condition on
$\Omega_\varepsilon \cap \{\Re z=0\}$ and keep the Steklov
condition on the
$\partial \Omega'_\varepsilon \cap \partial \Omega_\varepsilon.$
Let
$0=\sigma_0^N(\Omega'_\varepsilon) < \sigma_1^N(\Omega'_\varepsilon)
\le \sigma_2^N(\Omega'_\varepsilon) \dots $ be the eigenvalues of
this mixed problem (it is called  a {\it sloshing} problem, see
\cite{FK}).
Adding the Neumann condition inside the domain increases
the space of test functions,
and hence, by the standard monotonicity argument \cite[p. 100]{Bandle},
it pushes the eigenvalues down.
Therefore, $$\sigma_2(\Omega_\varepsilon) \ge \sigma_1^N(\Omega'_\varepsilon),$$ and
hence to prove \eqref{opposite} it suffices to show that
\begin{equation}
\label{slosh} \lim_{\varepsilon \to 0+}
\sigma_1^N(\Omega'_\varepsilon) =1.
\end{equation}
The family of domains $\Omega'_\varepsilon$ converges to $\D$ in the
Hausdorff
complementary topology  (see \cite[p. 101]{BB})
as
$\varepsilon \to 0+$. 
Moreover, since the domains $\Omega'_\varepsilon$ are uniformly Lipschitz in both $z\in \partial
\Omega'_\varepsilon$ and $\varepsilon\in (0, \frac{1}{10})$, the extension operators
$H^1(\Omega_\varepsilon) \to H^1(\mathbb{R}^2)$  are uniformly
bounded \cite[p. 198]{BB}, and the norms of the trace restriction
operators are uniformly bounded as well \cite{Di}. Note also that
the Neumann part of $\partial \Omega'_\varepsilon$ given by
$\Omega_\varepsilon \cap \{\Re z=0\}$ tends to the single point
$z=0$ as $\varepsilon \to 0+$. Therefore, using the Rayleigh
quotient for the sloshing problem \cite[p. 673]{FK}) we get
$$
\lim_{\varepsilon \to 0+}
\sigma_n^N(\Omega'_\varepsilon)=\sigma_n(\D),\,\, n=1,2,\dots
$$
in the same way as \cite[Corollary 7.4.2]{BB}. Taking $n=1$ we get
\eqref{slosh}. This completes the proof of the lemma.
\end{proof}

\smallskip

Let us now complete the proof of Theorem \ref{maintheorem}. First, it follows from \eqref{hps2} and the obvious inequality
$\sigma_{n+1}(\Omega_\varepsilon) \ge \sigma_n(\Omega_\varepsilon)$,
that \eqref{limitcase} implies \eqref{limitcasenew}. Therefore, it
suffices to prove \eqref{limitcase}. For $n=2$  it  follows from
Lemma \ref{conv}.  For $n>2$ the proof is analogous.  Define
$\Omega_\varepsilon$ as the union of $n$ disks of radius $1+\varepsilon$
centered at the points $z=2k$, $k=0,1,\dots,n-1$ (the domain
$\Omega_\varepsilon$ looks like a necklace).  We make cuts along the
vertical lines $\Re z=2k-1$, $k=1,2,\dots,n$, and impose Neumann boundary
conditions along these cuts.  We get  $n$ auxiliary mixed
problems. They are of two types:   the first and the last disks have
just one cut (we denote the corresponding domains by
$\Omega'_\varepsilon$ as before), and the intermediate disks have two
cuts --- one on the left and the other one on the right (the
corresponding domains are denoted by $\Omega''_\varepsilon)$.  The spectra of each
of the $n$ auxiliary problems start with the zero eigenvalue. Using
the same monotonicity  and convergence arguments  as above, we get
$$\sigma_n(\Omega_\varepsilon) \ge
\min\left(\sigma_1^N(\Omega'_\varepsilon),\, \sigma_1^N(\Omega''_\varepsilon)\right)$$
and
$$
 \lim_{\varepsilon \to 0+}
\sigma_1^N(\Omega'_\varepsilon) =  \lim_{\varepsilon \to 0+}
\sigma_1^N(\Omega''_\varepsilon) =1.
$$
Therefore,
$
\liminf_{\varepsilon \to 0+}
\sigma_n(\Omega_\varepsilon) \ge 1.
$
Since $\lim_{\varepsilon \to 0+} M(\Omega_\varepsilon) = 2\pi n$,  it follows from  \eqref{hps} that
$
\limsup_{\varepsilon \to 0+}
\sigma_n(\Omega_\varepsilon) \le 1.
$
Hence, $\lim_{\varepsilon \to 0+} \sigma_n(\Omega_\varepsilon)=1$ and
this completes the proof of Theorem~\ref{maintheorem}.

\subsection{Collapse of the Steklov spectrum: an example}
\label{counterexample}  One could ask why the sequence
$\Omega_\varepsilon$ is constructed by pulling the disks apart,
rather than joining them by a  tiny passage disappearing as
$\varepsilon \to 0$.  While this looks geometrically more natural, it
turns out that  the behavior of the Steklov spectrum under such
degeneration can be quite  unexpected.

As before, set $\rho \equiv 1$.
\begin{figure}[h]
  \centering
  \psfrag{e}[][][1]{$\varepsilon$}
  \psfrag{f}[][][1]{$\varepsilon^3$}
  \includegraphics[width=11cm]{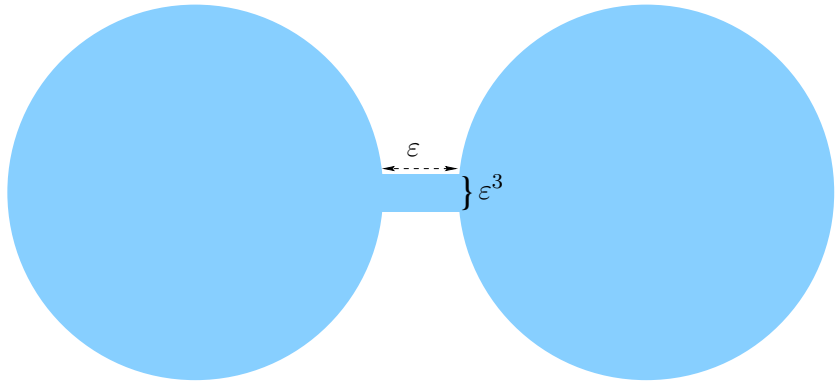}
  \caption{The domain $\Sigma_\varepsilon$}
\end{figure}
Let $\Sigma_\varepsilon=\D_1\cup
P_\varepsilon \cup \D_2$, where $\D_1$ and $\D_2$ are two
copies of the unit disk joined  by a rectangular passage
$P_\varepsilon$ of length $\varepsilon$ and width $\varepsilon^3$. What is essential in this construction 
is that the width of the passage tends to zero much faster than its length.
For simplicity we assume that the disks and the passage are chosen
in such a way that the domain
$\Sigma_\varepsilon$ is symmetric with respect to both coordinate
axes. Then, surprisingly enough,
\begin{equation}
\label{collapse} \lim_{\varepsilon \to 0+}
\sigma_n(\Sigma_\varepsilon)=0 \,\,\, {\rm for} \,\,\, {\it all}
\,\,\, n=1,2,\dots.
\end{equation}
Indeed, consider pairwise orthogonal test functions  vanishing in
$\D_1\cup \D_2$, and equal to $\sin \frac{2\pi nx}{\varepsilon}$ in
the passage $P_\varepsilon$.  For each $n$, the gradient of the test
function is of order $n/\varepsilon$, the area of $P_\varepsilon$ is
$\varepsilon^4$ and the length of the boundary of $P_\varepsilon$ is
$2\varepsilon$. Therefore, for each fixed $n$, the corresponding
Rayleigh quotient is of order $n^2 \varepsilon$  and tends to zero
as $\varepsilon \to 0+$.  This proves \eqref{collapse}.

Similar constructions were studied in the context of Neumann
boundary conditions (see \cite{JM, HSS} and references therein).
However, the Neumann eigenvalues of $\Sigma_\varepsilon$ converge to
the corresponding eigenvalues of the disjoint union of two disks as
$\varepsilon \to 0+$.  The total ``collapse'' of the Steklov
spectrum in the example above is caused by the fact that the
denominator of the Rayleigh quotient is an integral over the {\it
boundary}. Note that 
the perimeter of the passage $P_\varepsilon$ tends to zero much
slower than its area, and hence, for every fixed $n$, the numerator
in the Rayleigh quotient vanishes much faster than the denominator.

\bigskip

In the subsequent sections we prove Theorem~\ref{maintheoremII}.

\section{Folding and rearrangement of measure}
\label{folding:section}
\subsection{Conformal mapping to a disk}
\label{conform}
Let $\Omega$ be a simply connected planar domain with Lipschitz
boundary. As before,  $\mathbb{D}=\left\{z\in\mathbb{C}\, \bigl|\bigr.\,
|z|<1\right\}$ is the open unit disk. By the Riemann mapping theorem
(see \cite[p. 342]{Taylor}), there exists a conformal equivalence
$\phi:\mathbb{D}\rightarrow\Omega$ which extends to a  homeomorphism
$\overline{\mathbb{D}}\rightarrow\overline{\Omega}$  (slightly abusing notations, here and further on we  denote a
conformal map and its extension to the boundary by the same symbol).
Let $ds$ be the arc-length measure on $\partial\Omega$, and  $d\mu$
be the pull-back by $\phi$ of the measure $\rho(s)ds$:
\begin{gather}
\label{induced}
  \int_{\mathcal{O}}d\mu=\int_{\phi(\mathcal{O})}\rho(s)\,ds
\end{gather}
 for any open set $\mathcal{O}\subset S^1$.
Taking \eqref{induced} into account and using  conformal invariance
of the Dirichlet integral, we rewrite the variational
characterization~(\ref{rayleigh}) of $\sigma_2$ as follows:
\begin{gather}
\label{rayleighdisk}
 \sigma_2(\Omega)=\inf_{E} \sup_{0\neq u\in E}
  \frac{\int_{\mathbb{D}}|\nabla u|^2\,dz}
  {\int_{S^1}u^2\,d\mu}.
\end{gather}
Here the infimum is taken over all subspaces $E \subset
H^1(\mathbb{D})$, such that ${\rm dim}\,\, E=~2$ and
$\int_{S^1}u\,d\mu=0$ for all $u \in E$.

\subsection{Hyperbolic caps}
\label{hcaps} Let $\gamma$ be a geodesic in the Poincar\'e disk
model, that is a diameter or the intersection of the disk with a
circle which is orthogonal to $S^1$.
\begin{figure}[h]
  \centering
  \psfrag{a}[][][1]{$a_{l,p}$}
  \psfrag{p}[][][1]{$p$}
  \psfrag{l}[][][1]{$l$}
  \includegraphics[width=7cm]{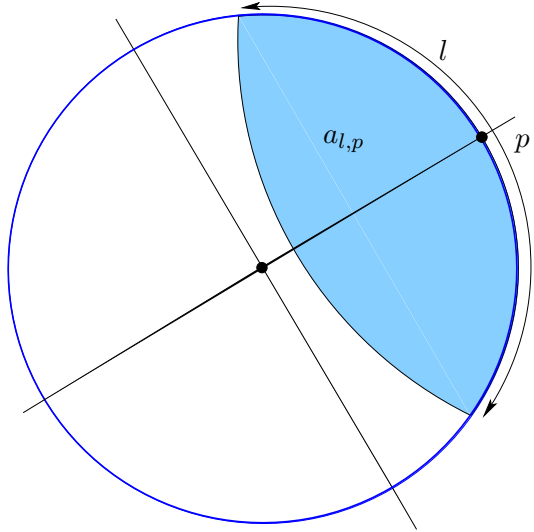}
  \caption{The hyperbolic cap $a_{l,p}$}
  \label{fighypocapo}
\end{figure}
Each connected component of
$\mathbb{D}\setminus\gamma$ is called a \emph{hyperbolic cap}
\cite{GNP}.
Given $p\in S^1$ and $l\in (0,2\pi)$, let $a_{l,p}$ be
the hyperbolic cap such that the circular segment $\partial
a_{l,p}\cap S^1$ has length $l$ and is centered at $p$
(see Figure~\ref{fighypocapo}). This gives
an identification of the space  $\mathcal{HC}$ of all hyperbolic
caps with the cylinder $(0,2\pi)\times S^1$. Given a cap $a \in
{\mathcal HC}$, let
$\tau_a:\overline{\mathbb{D}}\rightarrow\overline{\mathbb{D}}$ be
the reflection across the hyperbolic geodesic bounding $a$. That is,
$\tau_a$ is the unique non-trivial conformal involution of
$\mathbb{D}$ leaving every point of the geodesic $\partial
a\cap\mathbb{D}$ fixed. In particular,
$\tau_a(a)=\mathbb{D}\setminus\overline{a}.$ The {\it lift} of a
function $u:\overline{a}\rightarrow\mathbb{R}$ is the function
$\tilde{u}:\overline{\mathbb{D}}\rightarrow\mathbb{R}$ defined by
\begin{equation}
\label{lift} \tilde{u}(z)=
\begin{cases}
  u(z) & \mbox{if } z\in \overline{a},\\
  u(\tau_az) & \mbox{if } z\in\overline{\mathbb{D}\setminus a}.
\end{cases}
\end{equation}
Observe that
\begin{align}\label{folding}
  \int_{S^1}\tilde{u}\,d\mu&=\int_{\partial a\cap S^1}u\,d\mu+\int_{\tau_a(\partial a)\cap S^1}u\circ\tau_a\,d\mu\nonumber\\
  &=\int_{\partial a\cap S^1}u\,(d\mu+\tau_a^*d\mu).
\end{align}
The measure
\begin{equation}
\label{folded}
 d\mu_a =
\begin{cases}
  d\mu+\tau_a^*d\mu& \mbox{on } \partial a\cap S^1,\\
  0 & \mbox{on } S^1\setminus \partial a
\end{cases}
\end{equation}
is called the \emph{folded measure}. Equation~(\ref{folding}) can be rewritten
as
$$\int_{S^1}\tilde{u}\,d\mu=\int_{S^1}u\,d\mu_a.$$

\subsection{Eigenfunctions on the disk}
Given $t\in\mathbb{R}^2,$ define
$X_t:\overline{\mathbb{D}}\rightarrow\mathbb{R}$ by $X_t(z)=z\cdot t$, the
inner product of $z$ and $t$ in $\mathbb{R}^2$. 
Let $(e_1, e_2)$ be the standard basis of $\mathbb{R}^2$. Then
$X_{e_1}$ and $X_{e_2}$ form a basis of the first Steklov
eigenspace on the disk with $\rho \equiv 1$. Using Hersch
renormalization procedure (see \cite[subsection 4.1]{GNP}), we
assume that the center of mass of the measure $d\mu$ is at the
origin:
\begin{gather}\label{conventionI}
  \int_{S^1}X_t\,d\mu=0, \, \, \,   \forall t\in\mathbb{R}^2.
\end{gather}
Using a rotation if necessary, we may also assume that
\begin{gather}\label{conventionII}
    \int_{S^1}X_{e_1}^2\,d\mu\geq\int_{S^1}X_{t}^2 \,d\mu, \, \, \, \forall t\in S^1.
\end{gather}

\subsection{Rearranged measure}
Let $a\in\mathcal{HC}$ be a hyperbolic cap and let
$\psi_a:\mathbb{D}\rightarrow a$ be a conformal equivalence.
Following the convention adopted in subsection \ref{conform}, we also denote
its extension $\overline{\mathbb D} \to \overline{a}$ by
$\psi_a$.  For each $t\in\mathbb{R}^2$, define $u_a^t:\overline{a}
\rightarrow\mathbb{R}$ by
$$u_a^t(z)=X_t\circ\psi_a^{-1}(z)=t\cdot\psi_a^{-1}(z).$$
The following auxiliary lemma will be used in the proof of Lemma \ref{dirichlet}.
\begin{lemma}\label{notharmonic}
The lift of the function  ${u}_a^t$  is not harmonic in ${\mathbb D}$.
\end{lemma}
\begin{proof}
  Suppose that $\tilde{u}_a^t$ is harmonic. Then it is smooth, and by
  \eqref{lift}, the normal derivative of $u_a^t$ vanishes at any point
  $p\in\partial a\cap \mathbb{D}$. It is well-known that the vanishing of the normal derivative is
  preserved by conformal transformations. It follows that the normal derivative of the function
  $X_t=u_a^t\circ\psi_a$ vanishes on
  $\psi_a^{-1}(\partial a\cap\mathbb{D}) \subset S^1$.
  However, a straightforward computation shows that for any
  $s \neq \pm \frac{t}{|t|}$, $\frac{\partial}{\partial n}X_t(s) \neq 0.$
\end{proof}

Let $w_a^t\in C^{\infty}(\mathbb{D})$ be the unique
harmonic extension of $\tilde{u}_a^t\bigl|\bigr._{S^1}$, that is
\begin{gather}
\label{wat}
\begin{cases}
  \Delta w_a^t=0& \mbox{ in } \mathbb{D},\\
  w_a^t=\tilde{u}_a^t& \mbox{ on }S^1.
\end{cases}
\end{gather}
These functions will later be used as test functions in the variational
characterization~\ref{rayleighdisk}. Observe that
\begin{align}
\label{observe}
  \int_{S^1}\tilde{u}_a^t\,d\mu=\int_{S^1}u_a^t\,d\mu_a
  =\int_{S^1}X_t\,\psi_a^*d\mu_a.
\end{align}
We call the pullback measure
\begin{equation}
\label{rearranged}
d\nu_a=\psi_a^*d\mu_a
\end{equation}
 the \emph{rearranged measure} on $S^1$.

A family of conformal transformations $\left\{\psi_a:\mathbb{D}\rightarrow a\right\}_{a\in\mathcal{HC}}$ is said to be {\it continuous} if the map
$(0,2\pi)\times S^1\times\mathbb{D} \to {\mathbb D}$ defined by $(l,p,z)\mapsto\psi_{a_{l,p}}(z)$ is continuous. The next lemma describes the properties of the rearranged measure $d\nu_a$ as the cap $a$ degenerates either to the full disk or to a point $p\in S^1$.

\begin{lemma}\label{flipfloplemma}
  There exists a continuous family of conformal equivalences
  $\left\{\psi_a:\mathbb{D}\rightarrow a\right\}_{a\in\mathcal{HC}}$
  such that for each cap $a\in\mathcal{HC}$ and each
  $t\in\mathbb{R}^2$,
  \begin{gather}
    \int_{S^1}w_a^t\,d\mu=0,\label{normalization}\\
    \lim_{a\rightarrow \mathbb{D}}d\nu_a=d\mu,\label{fulldisk}\\
    \lim_{a\rightarrow p}d\nu_a=R_p^*d\mu\label{flipflopProperty},
  \end{gather}
  where $w_a^t$ is defined by \eqref{wat}, $d\nu_a$ is the rearranged measure given by \eqref{rearranged},
  $p \in S^1$ and $R_p(x)=x-2(x\cdot p)$ is the reflection with
respect to the diameter orthogonal to the vector  $p$.
\end{lemma}
Here and further on, the topology on measures is induced by the following norm:
\begin{equation}
\label{norm}
  \|d\nu\|=\sup_{f\in C(S^1), |f|\le 1} \left|\int_{S^1}f\,d\nu\,
  \right|.
\end{equation}

\begin{proof}
    Let us give an outline of the proof, for more details, see~\cite[Section 2.5]{GNP}.
    Start with any continuous family
    $\left\{\phi_a:\mathbb{D}\rightarrow a\right\}_{a\in\mathcal{HC}}$
    such that $\lim_{a\rightarrow\mathbb{D}}\phi_a=\id$.
    The maps $\psi_a$ are defined by composing the $\phi_a$'s on both sides with automorphisms of the disk
    appearing in the Hersch renormalization procedure.
    In particular,~\eqref{normalization} is automatically satisfied.
    As the cap $a$ converges to the full disk $\mathbb{D}$, the conformal equivalences $\psi_a$ converge
    to the identity map on $\mathbb{D}$, which implies~\eqref{fulldisk}.
    Finally, setting $n=1$ in~\cite[Lemma~4.3.2]{GNP} one gets \eqref{flipflopProperty}.
\end{proof}
From now on, we fix the family of conformal maps $\psi_a$ defined in
Lemma~\ref{flipfloplemma}.
Lemma \ref{flipfloplemma} implies that the rearranged measure $d\nu_a$ depends continously on the cap $a$.
This is essential  for the topological argument used  in the proof of Proposition \ref{multexist}.

\section{Construction of test functions}
\label{testfunctions}
\subsection{Estimate on the Rayleigh quotient} It follows from ~(\ref{normalization})
that the functions $w_a^t$ defined by \eqref{wat} are admissible in
the variational characterization~(\ref{rayleighdisk}) for
$\sigma_2$. For each hyperbolic cap $a \in {\mathcal HC}$ let
$$E_a=\left\{w_a^t\, \mid\,t\in\mathbb{R}^2\right\}.$$
be a two-dimensional space of test functions.
\begin{lemma}
\label{dirichlet} For any test function $w_a^t \in E_a$,
$$
 \int_{\mathbb{D}}|\nabla w_a^t|^2\,dz< 2\pi.
$$
\end{lemma}
\begin{proof}
It is well-known that a harmonic function, such as $w_a^t$, is the
unique minimizer of the Dirichlet energy among all functions with
the same boundary data~(see \cite[p. 157]{JST}). It follows from
Lemma~\ref{notharmonic} that $w_a^t\neq\tilde{u}_a^t$ in
$H^1(\mathbb{D})$. Therefore,
\begin{align}\label{doubling}
  \int_{\mathbb{D}}|\nabla w_a^t|^2\,dz&<
  \int_{\mathbb{D}}|\nabla\tilde{u}_a^t|^2\,dz
  =\int_{a}|\nabla u_a^t|^2\,dz+
  \int_{\mathbb{D}\setminus a}|\nabla (u_a^t\circ\tau_a)|^2\,dz\nonumber\\
  &=2\int_{a}|\nabla u_a^t|^2\,dz
  =2\int_{\mathbb{D}}|\nabla X_t|^2\,dz=2\underbrace{\sigma_1(\mathbb{D})}_1
  \overbrace{\int_{S^1}X_t^2\,d\theta}^{\pi}=2\pi,
\end{align}
where the second and the third equalities follow from the conformal
invariance of the Dirichlet energy.
\end{proof}

Let $t_1, t_2 \in S^1$ be such that $t_1 \cdot t_2 =0$. Given a
hyperbolic cap $a\in\mathcal{HC}$,  we have 
\begin{align}
\label{l2}
  \int_{S^1}(w_a^{t_1})^2\,d\mu&=\int_{S^1} (X_{t_1})^2 d\nu_a\nonumber\\
  &\geq
  \frac{1}{2}\int_{S^1}\overbrace{(X_{t_1})^2+(X_{t_2})^2}^1\,d\nu_a
  =\frac{1}{2}\int_{\partial\Omega}\rho(s)\,ds.
\end{align}
Here the first equality follows from \eqref{wat} and
\eqref{observe}, the last equality follows from \eqref{folded}
and \eqref{induced}, and we may assume without loss of generality
that the inequality in the middle is true (if not, we interchange $t_1$ and $t_2$).


\begin{remark}
Since $X_{t_1}^2+X_{t_2}^2=1$ on $S^1$, the estimate \eqref{l2}  is proved as in  \cite{Hersch},  and it is  much easier than the analogous result
\cite[Lemma 2.7.5]{GNP} for  the Neumann  problem.
\end{remark}

Consider the one-dimensional space of test functions
$$V_{t_1}=\left\{\alpha w_a^{t_1}\ \mid\ \alpha\in\mathbb{R}\right\}.$$
It follows from Lemma \ref{dirichlet} and formula \eqref{l2} that any function $u\in V_{t_1}$ satisfies
\begin{gather}
\label{rayl}
  \frac{\int_{\mathbb{D}}|\nabla u|^2\,dz}
  {\int_{S^1}u^2\,d\mu}\leq
  \frac{4\pi}{M(\Omega)}.
\end{gather}
Our next goal is to show that there exists a hyperbolic cap $a
\subset {\mathbb D}$ such that \eqref{rayl} holds not only for $u\in
V_{t_1}$ but for {\it each} $u\in E_a$. Since  $E_a$ is
two-dimensional, the estimate \eqref{bound:main} will follow
from~\eqref{rayl} and ~\eqref{rayleighdisk}.
\subsection{Simple and multiple measures}
Given a finite measure $d\nu$ on $S^1$, consider
the quadratic form  $V_{d\nu}:\mathbb{R}^2\rightarrow\mathbb{R}$ defined by
$$V_{d\nu}(t)=\int_{S^1}X_t^2\, d\nu.$$

Let $\mathbb{R}P^1=S^1/\mathbb{Z}_2$ be the projective line. We
denote by $[t]\in \mathbb{R}P^1$ the element of the projective line
corresponding to the pair of points $\pm t \in S^1$. We say that
$[t] \in \mathbb{R}P^1$ is a {\it maximizing direction} for the
measure $d\nu$ if $V_{d\nu}([t])\ge V_{d\nu}([s])$ for any $[s] \in \mathbb{R}P^1$. The measure
$d\nu$ is called {\it simple} if there is a unique maximizing
direction. Otherwise, the measure $d\nu$ is said to be {\it
multiple}.
\begin{lemma}
\label{multcaplemma}
  A measure $d\nu$ is multiple if and only if $V_{d\nu}(t)$ does not depend on
  $t \in S^1$.
\end{lemma}
\begin{proof} The lemma follows from the fact that $V_{d\nu}(t)$ is a
quadratic form, and is proved  analogously to \cite[Lemma
2.6.1]{GNP}.
\end{proof}
Note that by \eqref{conventionII}, $[e_1]$ is a maximizing
direction for the measure $d\mu$.
\begin{proposition} \label{multexist} If the measure $d\mu$ is simple, then there exists
  a cap $a\in \mathcal{HC}$  such that the rearranged measure $d\nu_a$ is
  multiple.
\end{proposition}

Proposition \ref{multexist} is proved  by contradiction. Assume that the measure $d\mu$, as
well as the measures $d\nu_a$ for all $a\in {\mathcal HC}$, are simple. Given a hyperbolic cap
$a$,  let $[m(a)]\in\mathbb{R}P^1$ be the unique
maximizing direction for $d\nu_a$.

By construction, the folded measures $d\mu_a$ depend continuously
on the cap $a$. The family $\psi_a$ is continuous by Lemma \ref{flipfloplemma}, and hence
the rearranged measures $d\nu_a$ depend continuously on
$a$. Therefore, the functions
$V_{d\nu_a}$ and  the unique maximizing
direction $[m(a)]$ also depend continuously on $a$.

Let us understand the behavior of the maximizing
direction as the cap $a$ degenerates either to the full disk or to a
point.
\begin{lemma}\label{degeneratecaps}
  Let  the measure $d\mu$ as well as the measures
  $d\nu_a$ for all $a\in {\mathcal HC}$  be simple.
  Then
  \begin{gather}
    \lim_{a\rightarrow\mathbb{D}} [m(a)]=[e_1]\label{eqndegenerate1}\\
    \lim_{a\rightarrow e^{i\theta}} [m(a)]=
    [e^{2i\theta}]\label{eqndegenerate2}.
  \end{gather}
\end{lemma}
\begin{proof}
  First, note that formula~\eqref{eqndegenerate1} immediately follows from
  \eqref{fulldisk} and \eqref{conventionII}.
  Let us prove \eqref{eqndegenerate2}. Set $p=e^{i\theta}$.
  Formula~\eqref{flipflopProperty} implies
  \begin{equation}
    \label{Xps}
    \lim_{a\rightarrow p}\int_{S^1}X_t^2\,d\nu_a=
    \int_{S^1}X_t^2\,R_p^*d\mu
    =\int_{S^1}X_t^2\circ R_p\,d\mu
    =\int_{S^1}X_{R_pt}^2\,d\mu.
  \end{equation}
  Since $d\mu$ is simple, $[e_1]$ is the unique maximizing direction for $d\mu$
  and the right hand side of \eqref{Xps} is maximal for
  $R_pt=\pm e_1$. Applying $R_p$ on both sides we get
  $t=\pm e^{2i\theta}$ and hence  $[m(a)]=[e^{2i\theta}]$.
\end{proof}
\begin{proof}[Proof of Proposition~\ref{multexist}]
  Suppose that for each hyperbolic cap $a\in\mathcal{HC}$ the measure $d\nu_a$
  is simple.
  Recall that the space $\mathcal{HC}$ is identified with the open cylinder
  $(0,2\pi)\times S^1$.
  Define
  $h:(0,2\pi)\times S^1\rightarrow \mathbb{R}P^1$
  by $h(l,p)=[m(a_{l,p})].$ As was mentioned above, the maximizing direction depends continuously
  on the cap $a$.
  Therefore, it follows from Lemma~\ref{degeneratecaps}  that  $h$ extends to a continuous
  map on the closed cylinder $[0,2\pi]\times S^1$ such that
  $$h(0,e^{i\theta})=[e_1],\ h(2\pi,e^{i\theta})=[e^{2i\theta}].$$
  This means that $h$ is a homotopy between a trivial loop and a
  non-contractible loop on $\mathbb{R}P^1$. This is a contradiction.
\end{proof}

\subsection{Proof of Theorem~\ref{maintheoremII}}
  Assume that the measure $d\mu$ is simple. By
  Proposition~\ref{multexist},  there exists a cap
  $a\in\mathcal{HC}$ such that the measure $d\nu_a$ is multiple so that
  inequality~\eqref{rayl} holds for any $u\in E_a$.
  Theorem~\ref{maintheoremII} then immediately follows from the variational characterization~(\ref{rayleighdisk}) of
  $\sigma_2$.

  Suppose now that the measure $d\mu$ is multiple. In this case the proof is easier.
  Indeed, it follows from Lemma~\ref{multcaplemma}, that
  any direction $[s]\in\mathbb{R}P^1$ is maximizing for $d\mu$ so that
  we can use the space
  $$E=\left\{X_t \,\bigl|\bigr. t\in\mathbb{R}^2\right\}$$
  of test functions in the variational
  characterization~(\ref{rayleighdisk}) of $\sigma_2.$
  Replacing $w_a^t$ by $X_t$ and inspecting~\eqref{doubling} we notice
  that the factor $2$ disappears. Therefore, \eqref{rayleighdisk}
  implies
  \begin{equation}
    \label{better}
    \sigma_2(\Omega)\, M(\Omega) \le  2\pi,
  \end{equation}
  which is an even better bound than~\eqref{bound:main}. This
  completes the proof of Theorem \ref{maintheoremII}.
\begin{remark}
  When $d\mu$ is multiple, Lemma~\ref{notharmonic} is not applicable, since we
  are not using formula \eqref{lift}, and therefore the inequality
  \eqref{better} is not strict. Indeed, the equality is attained on
  a disk with $\rho \equiv {\rm const}$.

It is easy to show that if the domain $\Omega$ is
symmetric of order $q\ge 3$ in the sense of  \cite{Bandle2} and \cite[pp.
136-140]{Bandle}  (for instance, if $\Omega$ is a regular
$q$-gon), then the measure $d\mu$ is multiple, provided the density
$\rho$ satisfies the same symmetry condition. Under these
assumptions \eqref{better} is a special case of \cite[Theorem
3.15]{Bandle}. In fact, one can show using Courant's nodal domain
theorem for Steklov eigenfunctions \cite[section 3]{KS} that if the
domain $\Omega$ and the density $\rho$ are symmetric of order $q$,
then $\sigma_1=\sigma_2$ so that \eqref{better} is just a
consequence of \eqref{Weinst}. Indeed, in this case $\Omega$ has at
least two axes of symmetry, and each of them is a nodal line of an
eigenfunction corresponding to $\sigma_1$. Therefore, ${\rm
mult}(\sigma_1)\ge 2$. We are not aware of any examples for which
\eqref{better} gives new information, i.e. the measure $d\mu$ is
multiple but $\sigma_2  > \sigma_1$.
\end{remark}


\subsection*{Acknowledgments}  We are thankful to Marl\`ene Frigon,
Michael Levitin, Marco Marletta, Nikolai Nadirashvili  and Yuri
Safarov for useful discussions.

\end{document}